\documentclass[a4paper,11pt]{amsart}
\usepackage{graphicx}
\usepackage{mathptmx}
\usepackage{amsmath}
\usepackage{amssymb}
\usepackage{enumitem}
\usepackage{xcolor}
\usepackage{xparse}
\NewDocumentCommand{\eulerian}{omm}
 {%
  \genfrac<>{0pt}{}{#2}{#3}%
  \IfValueT{#1}{_{\!#1}}%
 }

 \newmuskip\pFqmuskip

\newcommand*\pFq[6][8]{%
  \begingroup 
  \pFqmuskip=#1mu\relax
  \mathchardef\normalcomma=\mathcode`,
  \mathcode`\,=\string"8000
  \begingroup\lccode`\~=`\,
  \lowercase{\endgroup\let~}\pFqcomma
  {}_{#2}F_{#3}{\left(\genfrac..{0pt}{}{#4}{#5}\bigg|#6\right)}%
  \endgroup
}
\newcommand{\pFqcomma}{{\normalcomma}\mskip\pFqmuskip}

\newtheorem{theorem}{Theorem}

\begin{document}

\title[On generalized degenerate Euler-Genocchi polynomials]{On generalized degenerate Euler-Genocchi polynomials}

\author{Taekyun  Kim*}
\address{Department of Mathematics, Kwangwoon University, Seoul 139-701, Republic of Korea}
\email{tkkim@kw.ac.kr}

\author{DAE SAN KIM*}
\address{Department of Mathematics, Sogang University, Seoul 121-742, Republic of Korea}
\email{dskim@sogang.ac.kr}

\author{Hye Kyung Kim*}
\address{Department of Mathematics Education, Daegu Catholic University, Gyeongsan 38430,
Republic of Kore}
\email{hkkim@cu.ac.kr}

\subjclass[2010]{11B68; 11B73; 11B83}
\keywords{generalized degenerate Euler-Genocchi polynomials; generalized degenerate Euler-Genocchi polynomials of order $\alpha$; the alternating degenerate power sum of integers}
\thanks{ *corresponding authors}

\begin{abstract}
We introduce the generalized degenerate Euler-Genocchi polynomials as a degenerate version of the Euler-Genocchi polynomials. In addition, we introduce their higher-order version, namely the generalized degenerate Euler-Genocchi polynomials of order $\alpha$,  as a degenerate version of the generalized Euler-Genocchi polynomials of order $\alpha$. The aim of this paper is to study certain properties and identities involving those polynomials, the generalized falling factorials, the degenerate Euler polynomials of order $\alpha$, the degenerate Stirling numbers of the second kind, and the `alternating degenerate power sum of integers'.
\end{abstract}

\maketitle

\section{Introduction}
Explorations for various degenerate versions of some special numbers and polynomials have drawn the attention of many mathematicians in recent years, which were initiated by Carlitz's work in [4,5]. Many interesting results were obtained by exploiting different tools such as generating functions, combinatorial methods, $p$-adic analysis, umbral calculus, differential equations, probability theory, operator theory, analytic number theory and quantum physics (see [8-10,12-14,16-21]). \par
Belbachir et al. introduced the Euler-Genocchi polynomials in [2,3] and Goubi generalized them to the generalized Euler-Genocchi polynomials of order $\alpha$ in [7]. Here we introduce degenerate versions for both of them. Namely, we introduce the generalized degenerate Euler-Genocchi polynomials as a degenerate version of the Euler-Genocchi polynomials. In addition, we introduce their higher-order version, namely the generalized degenerate Euler-Genocchi polynomials of order $\alpha$, as a degenerate version of the generalized Euler-Genocchi polynomials of order $\alpha$. The aim of this paper is to study certain properties and identities involving those polynomials, the generalized falling factorials, the degenerate Euler polynomials of order $\alpha$, the degenerate Stirling numbers of the second kind, and the `alternating degenerate power sum of integers' (see \eqref{31}). The novelty of this paper is that it is the first paper that introduces the generalized degenerate Euler-Genocchi polynomials and the generalized degenerate Euler-Genocchi polynomials of order $\alpha$, as degenerate versions of the polynomials introduced earlier in [2,3,7]. \par
The outline of this paper is as follows. In Section1, we recall the degenerate exponentials, the degenerate Euler polynomials, the degenerate Euler polynomials of order $\alpha$, and the degenerate Genocchi polynomials of order $\alpha$. Also, we remind the reader of the degenerate Stirling numbers of the second and the incomplete Bell polynomials. Section 2 is the main result of this paper. We introduce the generalized degenerate Euler-Genocchi polynomials $A_{n,\lambda}^{(r)}(x)$, as a generalization of both the degenerate Euler polynomials and the degenerate Genocchi polynomials. In Theorem 1, the generalized falling factorials $(x)_{n,\lambda}$ are expressed in terms of $A_{n,\lambda}^{(r)}(x)$. A distribution property is derived for $A_{n,\lambda}^{(r)}(x)$ in Theorem 2. Then the generalized degenerate Euler-Genocchi polynomials of order $\alpha$, $A_{n,\lambda}^{(r,\alpha)}(x)$, are introduced as a higher-order version of  $A_{n,\lambda}^{(r)}(x)$. We deduce a simple expression for $A_{n,\lambda}^{(r,-m)}(x)$, with $m$ a positive integer in Theorem 3. In Theorems 4 and 5, we observe certain relations between $A_{n,\lambda}^{(r,\alpha)}(x)$ and the degenerate Euler polynomials of order $\alpha$, $\mathcal{E}_{n,\lambda}^{(\alpha)}(x)$. In Theorem 6, $\mathcal{E}_{n,\lambda}^{(\alpha)}(x)$ are expressed in terms of the degenerate Stirling numbers of the second kind $S_{2,\lambda}(n,k)$. In Theorem 7, $(x)_{n,\lambda}$ are represented in terms of $A_{n,\lambda}^{(r,\alpha)}(x)$  and $S_{2,\lambda}(n,k)$. In Theorem 8, we obtain an identity involving $A_{n,\lambda}^{(r,\alpha)}(x)$, $A_{n,\lambda}^{(r)}(x)$, and $\mathcal{E}_{n,\lambda}^{(\alpha-1)}(x)$ . Let $T_{k,\lambda}(n)$ denote the `alternating degenerate power sum of integers.' In Theorem 9, we get a representation of $T_{k,\lambda}(n)$ in terms of $S_{2,\lambda}(n,k)$. Finally, we get an expression of $A_{n,\lambda}^{(r)}(x)$ in terms of $T_{k,\lambda}(n)$ and $S_{2,\lambda}(n,k)$ in Theorem 10. In the rest of this section, we recall the facts that are needed throughout this paper. \par

\vspace{0.2cm}

It is well known that the Euler polynomials are defined by
\begin{equation}
	\frac{2}{e^{t}+1}e^{xt}=\sum_{n=0}^{\infty}E_{n}(x)\frac{t^{n}}{n!},\quad (\mathrm{see}\ [1-22]). \label{1}
\end{equation}
The Genocchi polynomials are given by
\begin{equation}
\frac{2t}{e^{t}+1}e^{xt}=\sum_{n=0}^{\infty}G_{n}(x)\frac{t^{n}}{n!},\quad (\mathrm{see}\ [1-22]),\label{2}
\end{equation}
When $x=0$, $G_{n}=G_{n}(0)$ are called the Genocchi numbers. \par
For any nonzero $\lambda\in\mathbb{R}$, the degenerate exponentials are defined by
\begin{equation*}
e_{\lambda}^{x}(t)=(1+\lambda t)^{\frac{x}{\lambda}}=\sum_{n=0}^{\infty}(x)_{n,\lambda}\frac{t^{n}}{n!},\quad e_{\lambda}(t)=e_{\lambda}^{1}(t),\quad (\mathrm{see}\ [12]),	
\end{equation*}
where the generalized falling factorials are given by
\begin{equation*}
(x)_{0,\lambda}=1,\quad (x)_{n,\lambda}=x(x-\lambda)\cdots(x-(n-1)\lambda),\quad (n\ge 1).
\end{equation*} \par
In [4,5], Carlitz introduced the degenerate Euler polynomials which are given by
\begin{equation*}
\frac{2}{e_{\lambda}(t)+1}e_{\lambda}^{x}(t)=\sum_{n=0}^{\infty}\mathcal{E}_{n,\lambda}(x)\frac{t^{n}}{n!}.
\end{equation*}
When $x=0$, $\mathcal{E}_{n,\lambda}=\mathcal{E}_{n,\lambda}(0)$ are called the degenerate Euler numbers. \par
For any nonzero $\alpha\in\mathbb{C}$, the degenerate Euler polynomials of order $\alpha$ are defined by
\begin{equation*}
\Big(\frac{2}{e_{\lambda}(t)+1}\Big)^{\alpha}e_{\lambda}^{x}(t)=\sum_{n=0}^{\infty}\mathcal{E}_{n,\lambda}^{(\alpha)}(x)\frac{t^{n}}{n!},\quad (\mathrm{see}\ [4,5]).
\end{equation*}
When $x=0$, $\mathcal{E}_{n,\lambda}^{(\alpha)}=\mathcal{E}_{n,\lambda}^{(\alpha)}(0)$ are called the degenerate Euler numbers of order $\alpha$. \par
Recently, the degenerate Genocchi polynomials of order $\alpha$ are defined by
\begin{equation*}
\Big(\frac{2t}{e_{\lambda}(t)+1}\Big)^{\alpha}e_{\lambda}^{x}(t)=\sum_{n=0}^{\infty}G_{n,\lambda}^{(\alpha)}(x)\frac{t^{n}}{n!},\quad (\mathrm{see}\ [18-22]).
\end{equation*}
When $x=0$, $G_{n,\lambda}^{(\alpha)}=G_{n,\lambda}^{(\alpha)}(0)$ are called the degenerate Genocchi numbers of order $\alpha$. \\
In particular, $\alpha=1$, $G_{n,\lambda}(x)=G_{n,\lambda}^{(1)}(x)$ are called the degenerate Genocchi polynomials. \par
For $n\ge 0$, the degenerate Stirling numbers of the second kind are introduced by Kim-Kim as
\begin{equation*}
(x)_{n,\lambda}=\sum_{k=0}^{n}S_{2,\lambda}(n,k)(x)_{k},\quad (\mathrm{see}\ [9]),	
\end{equation*}
where
\begin{equation*}
(x)_{0}=1,\quad (x)_{n}=x(x-1)\cdots(x-n+1),\quad (n\ge 1).\end{equation*} \par
For $k\ge 0$, the incomplete Bell polynomials are defined by
\begin{equation*}
\frac{1}{k!}\bigg(\sum_{i=1}^{\infty}x_{i}\frac{t^{i}}{i}\bigg)^{k}=\sum_{n=k}^{\infty}B_{n,k}(x_{1},x_{2},\dots,x_{n-k+1})\frac{t^{n}}{n!},\quad (\mathrm{see}\ [6,14,15,17 ]), 	
\end{equation*}
where
\begin{equation*}
B_{n,k}(x_{1},x_{2},\dots,x_{n-k+1})=\sum_{\substack{l_{1}+\cdots+l_{n-k+1}=k\\ l_{1}+2l_{2}+\cdots+(n-k+1)l_{n-k+1}=n}}\frac{n!}{l_{1}!l_{2}!\cdots l_{n-k+1}!}\bigg(\frac{x_{1}}{1!}\bigg)^{l_{1}}\cdots \bigg(\frac{x_{n-k+1}}{(n-k+1)!}\bigg)^{l_{n-k+1}}.
\end{equation*}

\section{Generalized degenerate Euler-Genocchi numbers and polynomials}
For $r\in\mathbb{Z}$ with $r\ge 0$, we consider the {\it{generalized degenerate Euler-Genocchi polynomials}} given by
\begin{equation}
\frac{2t^{r}}{e_{\lambda}(t)+1}e_{\lambda}^{x}(t)=\sum_{n=0}^{\infty}A_{n,\lambda}^{(r)}(x)\frac{t^{n}}{n!}. \label{3}	
\end{equation}
Note that $A_{0,\lambda}^{(r)}=A_{1,\lambda}^{(r)}(x)=\cdots=A_{r-1,\lambda}^{(r)}(x)=0$. \\
When $x=0$, $A_{n,\lambda}^{(r)}=A_{n,\lambda}^{(r)}(0)$ are called the {\it{generalized degenerate Euler-Genocchi numbers}}. Observe that
\begin{equation}
A_{n,\lambda}^{(0)}(x)=\mathcal{E}_{n,\lambda}(x),\quad A_{n,\lambda}^{(1)}(x)=G_{n,\lambda}(x),\quad (n\ge 0).\label{4}	
\end{equation} \par
From \eqref{3}, we have
\begin{equation}
A_{n,\lambda}^{(r)}(x+1)=\sum_{l=0}^{n}\binom{n}{l}(1)_{n-l,\lambda}A_{l,\lambda}^{(r)}(x),\quad (n\ge 0).
\end{equation}
By \eqref{3}, we get
\begin{align}
\sum_{n=0}^{\infty}(x)_{n,\lambda}\frac{t^{n}}{n!}&=e_{\lambda}^{x}(t)=\frac{1}{2t^{r}}\sum_{l=0}^{\infty}A_{l,\lambda}^{(r)}(x)\frac{t^{l}}{l!}\big(e_{\lambda}(t)+1\big) \label{6} \\
&=\frac{1}{2}\sum_{l=0}^{\infty}A_{l+r,\lambda}^{(r)}(x)\frac{t^{l}}{(l+r)!}\bigg(\sum_{m=0}^{\infty}\frac{(1)_{m,\lambda}}{m!}t^{m}+1\bigg)\nonumber \\
&=\frac{1}{2}\sum_{n=0}^{\infty}\sum_{l=0}^{n}\binom{n+r}{l+r}\frac{n!}{(n+r)!}(1)_{n-l,\lambda}A_{l+r,\lambda}^{(r)}(x)\frac{t^{n}}{n!}+\frac{1}{2}\sum_{n=0}^{\infty}A_{n+r,\lambda}^{(r)}(x)\frac{n!}{(n+r)!}\frac{t^{n}}{n!} \nonumber \\
&=\sum_{n=0}^{\infty}\frac{1}{2}\bigg(\sum_{l=0}^{n}\binom{n+r}{l+r}(1)_{n-l,\lambda}\frac{A_{l+r}^{(r)}(x)}{(n+r)_{r}}+\frac{A_{n+r,\lambda}^{(r)}(x)}{(n+r)_{r}}\bigg)\frac{t^{n}}{n!}.\nonumber	
\end{align}
Therefore, by comparing the coefficients on both sides of \eqref{6}, we obtain the following theorem.
\begin{theorem}
For $n\ge 0$, we have
\begin{equation*}
(x)_{n,\lambda}=\frac{1}{2(n+r)_{r}}\bigg(\sum_{l=0}^{n}\binom{n+r}{l+r}(1)_{n-l,\lambda}A_{l+r,\lambda}^{(r)}(x)+A_{n+r,\lambda}^{(r)}(x)\bigg).
\end{equation*}
\end{theorem}
For $m\in\mathbb{N}$ with $m\equiv 1\ (\mathrm{mod}\ 2)$, we have
\begin{align}
\sum_{n=0}^{\infty}A_{n,\lambda}^{(r)}(x)\frac{t^{n}}{n!}&=\frac{2t^{r}}{e_{\lambda}(t)+1}e_{\lambda}^{x}(t)=\frac{2t^{r}}{e_{\lambda}^{m}(t)+1}\sum_{l=0}^{m-1}	(-1)^{l}e_{\lambda}^{l+x}(t)\label{7} \\
&=\frac{2(mt)^{r}}{(e_{\lambda/m}(mt)+1)}\frac{1}{m^{r}} \sum_{l=0}^{m-1}(-1)^{l}e_{\lambda/m}^{\frac{l+x}{m}}(mt)\nonumber \\
&=\frac{1}{m^{r}}\sum_{l=0}^{m-1}(-1)^{l}\frac{2(mt)^{r}}{e_{\lambda/m}(mt)+1}e_{\lambda/m}^{\frac{l+x}{m}}(mt)\nonumber \\
&=\frac{1}{m^{r}}\sum_{l=0}^{m-1}(-1)^{l}\sum_{n=0}^{\infty}A_{n,\lambda/m}^{(r)}\Big(\frac{l+x}{m}\Big)m^{n}\frac{t^{n}}{n!}\nonumber \\
&=\sum_{n=0}^{\infty}m^{n-r}\sum_{l=0}^{m-1}(-1)^{l}A_{n,\lambda/m}^{(r)}\Big(\frac{l+x}{m}\Big)\frac{t^{n}}{n!}. \nonumber
\end{align}
Therefore, by comparing the coefficients on both sides of \eqref{7}, we obtain the following theorem.
\begin{theorem}
For $m\in\mathbb{N}$ with $m\equiv 1\ (\mathrm{mod}\ 2)$, we have
\begin{displaymath}
A_{n,\lambda}^{(r)}(x)= m^{n-r}\sum_{l=0}^{m-1}(-1)^{l}A_{n,\lambda/m}^{(r)}\Big(\frac{l+x}{m}\Big).
\end{displaymath}
\end{theorem}
For nonzero $\alpha \in\mathbb{C}$, and $r\in\mathbb{Z}$ with $r\ge 0$, we consider the {\it{generalized degenerate Euler-Genocchi polynomials of order $\alpha$}} which are given by
\begin{equation}
t^{r}\Big(\frac{2}{e_{\lambda}(t)+1}\Big)^{\alpha}e_{\lambda}^{x}(t)=\sum_{n=0}^{\infty}A_{n,\lambda}^{(r,\alpha)}(x)\frac{t^{n}}{n!}.\label{8}
\end{equation}
When $x=0$, $A_{n,\lambda}^{(r,\alpha)}=A_{n,\lambda}^{(r,\alpha)}(0)$ are called the {\it{generalized degenerate Euler-Genocchi numbers of order $\alpha$}}.  \par
From \eqref{8}, we note that
\begin{align}
A_{n,\lambda}^{(r,\alpha)}(x)&=\sum_{k=0}^{n}\binom{n}{k}A_{k,\lambda}^{(r,\alpha)}(x)_{n-k,\lambda}\label{10}\\
&=\sum_{k=0}^{n}\binom{n}{k}A_{n-k,\lambda}^{(r,\alpha)}(x)_{k,\lambda},\quad (n\ge 0). \nonumber
\end{align}
Let $\alpha=-m\ (m\in\mathbb{N})$. Then, by \eqref{8}, we get
\begin{align}
\sum_{n=0}^{\infty}A_{n,\lambda}^{(r,-m)}(x)\frac{t^{n}}{n!}&
=\frac{t^{r}}{2^{m}}(e_{\lambda}(t)+1)^{m}e_{\lambda}^{x}(t)=\frac{t^{r}}{2^{m}}\sum_{k=0}^{m}\binom{m}{k}e_{\lambda}^{k+x}(t) \label{11}\\	
&=\sum_{n=0}^{\infty}\frac{1}{2^{m}}\sum_{k=0}^{m}\binom{m}{k}(k+x)_{n,\lambda}\frac{t^{n+r}}{n!}\nonumber \\
&=\sum_{n=r}^{\infty}\frac{1}{2^{m}}\sum_{k=0}^{m}\binom{m}{k}(k+x)_{n-r,\lambda}(n)_{r}\frac{t^{n}}{n!}. \nonumber
\end{align}
Therefore, by comparing the coefficients on both sides of \eqref{11}, we obtain the following theorem.
\begin{theorem}
For $m\in\mathbb{N}$, we have
\begin{displaymath}
A_{n,\lambda}^{(r,-m)}(x)=\frac{(n)_{r}}{2^{m}}\sum_{k=0}^{m}\binom{m}{k}(k+x)_{n-r,\lambda}.
\end{displaymath}\par
When $x=0$, we have
\begin{equation}
A_{n,\lambda}^{(r,-m)}=\frac{(n)_{r}}{2^{m}}\sum_{k=0}^{m}\binom{m}{k}(k)_{n-r,\lambda}. \label{12}
\end{equation}
\end{theorem}
\noindent From \eqref{10} and \eqref{12}, we have
\begin{equation}
\begin{aligned}
A_{n,\lambda}^{(r,-m)}(x)&=\sum_{k=0}^{n}\binom{n}{k}A_{n-k,\lambda}^{(r,-m)}(x)_{k,\lambda}\\
&=\frac{1}{2^{m}}\sum_{k=0}^{n}\sum_{j=0}^{m}\binom{n}{k}\binom{m}{j}(j)_{n-k-r,\lambda}(n-k)_{r}(x)_{k,\lambda}.
\end{aligned}	\label{13}
\end{equation} \par
By \eqref{8}, we get
\begin{align}
\sum_{n=r}^{\infty}A_{n,\lambda}^{(r,\alpha)}(x)\frac{t^{n}}{n!}&=t^{r}\Big(\frac{2}{e_{\lambda}(t)+1}\Big)^{\alpha}e_{\lambda}^{x}(t) \label{14}\\
&=t^{r}\sum_{n=0}^{\infty}\mathcal{E}_{n,\lambda}^{(\alpha)}(x)\frac{t^{n}}{n!}=\sum_{n=r}^{\infty}\mathcal{E}_{n-r,\lambda}^{(\alpha)}(x)\frac{n!}{(n-r)!}\frac{t^{n}}{n!} \nonumber \\
&=\sum_{n=r}^{\infty}(n)_{r}\mathcal{E}_{n-r,\lambda}^{(\alpha)}(x)\frac{t^{n}}{n!}. \nonumber	
\end{align}
\begin{theorem}
For $n,r\ge 0$ with $n\ge r$, we have
\begin{displaymath}
A_{n,\lambda}^{(r,\alpha)}(x)=(n)_{r}\mathcal{E}_{n-r,\lambda}^{(\alpha)}(x).
\end{displaymath}
In particular, for $x=0$, we get
\begin{equation}\label{15}
A_{n,\lambda}^{(r,\alpha)}=(n)_{r}\mathcal{E}_{n-r,\lambda}^{(\alpha)}. 	
\end{equation}
\end{theorem}
By \eqref{10} and \eqref{15}, we get
\begin{align}
A_{n,\lambda}^{(r,\alpha)}(x)&=\sum_{k=0}^{n}\binom{n}{k}A_{n-k,\lambda}^{(r,\alpha)}(x)_{k,\lambda}=\sum_{k=0}^{n-r}\binom{n}{k}A_{n-k,\lambda}^{(r,\alpha)}(x)_{k,\lambda}
\label{16}	\\
&=\sum_{k=0}^{n-r}\binom{n}{k}(n-k)_{r}\mathcal{E}_{n-k-r,\lambda}^{(\alpha)}(x)_{k,\lambda}=(n)_{r}\sum_{k=0}^{n-r}\binom{n-r}{k}\mathcal{E}_{n-k-r,\lambda}^{(\alpha)}(x)_{k,\lambda}\nonumber \\
&=(n)_{r}(x)_{n-r,\lambda}+(n)_{r}\sum_{k=0}^{n-r-1}\binom{n-r}{k}\mathcal{E}_{n-k-r,\lambda}^{(\alpha)}(x)_{k,\lambda}.
\end{align}
Therefore, by \eqref{16}, we obtain the following theorem.
\begin{theorem}
For any nonzero $\alpha \in\mathbb{C}$ and $n,r\ge 0$ with $n\ge r$, we have
\begin{displaymath}
A_{n,\lambda}^{(r,\alpha)}(x)=(n)_{r}(x)
_{n-r,\lambda}+(n)_{r}\sum_{k=0}^{n-r-1}\binom{n-r}{k}\mathcal{E}_{n-k-r,\lambda}^{(\alpha)}(x)_{k,\lambda}.
\end{displaymath}
\end{theorem}
Let $f(t)=\sum_{n=0}^{\infty}a_{n}\frac{t^{n}}{n!}\in\mathbb{C}[\![t]\!]$, where $a_{0} \ne 0$.\\
For any nonzero $\alpha \in\mathbb{C}$, we have
\begin{align}
f^{\alpha}(t)&=\bigg(\sum_{k=0}^{\infty}a_{k}\frac{t^{k}}{k!}\bigg)^{\alpha}=\bigg(a_{0}+\sum_{k=1}^{\infty}a_{k}\frac{t^{k}}{k!}\bigg)^{\alpha} \label{17}\\
&=\sum_{k=0}^{\infty}\binom{\alpha}{k}a_{0}^{\alpha-k}\bigg(\sum_{i=1}^{\infty}a_{i}\frac{t^{i}}{i!}\bigg)^{k}\nonumber \\
&=a_{0}^{\alpha}+\sum_{k=1}^{\infty}(\alpha)_{k}a_{0}^{\alpha-k}\frac{1}{k!}\bigg(\sum_{i=1}^{\infty}a_{i}\frac{t^{i}}{i!}\bigg)^{k}. \nonumber\\
&=a_{0}^{\alpha}+\sum_{k=1}^{\infty}(\alpha)_{k}a_{0}^{\alpha-k}\sum_{n=k}^{\infty}B_{n,k}(a_{1},a_{2},\dots,a_{n-k+1})\frac{t^{n}}{n!} \nonumber \\
&=a_{0}^{\alpha}+\sum_{n=1}^{\infty}\sum_{k=1}^{n}(\alpha)_{k}a_{0}^{\alpha-k}B_{n,k}(a_{1},\dots,a_{n-k+1})\frac{t^{n}}{n!}.\nonumber
\end{align} \par
From the definition of degenerate Stirling numbers of the second kind, we have
\begin{align}
\sum_{n=k}^{\infty}S_{2,\lambda}(n,k)\frac{t^{n}}{n!}&=\frac{1}{k!}(e_{\lambda}(t)-1)^{k} \label{18}\\
&=\frac{1}{k!}\sum_{j=0}^{k}\binom{k}{j}(-1)^{k-j}\sum_{n=0}^{\infty}(j)_{n,\lambda}\frac{t^{n}}{n!}\nonumber \\
&=\sum_{n=0}^{\infty}\frac{1}{k!}\sum_{j=0}^{k}\binom{k}{j}(-1)^{k-j}(j)_{n,\lambda}\frac{t^{n}}{n!}.\nonumber
\end{align}
Comparing the coefficients on both sides of \eqref{18}, we have
\begin{displaymath}
S_{2,\lambda}(n,k)=\frac{1}{k!}\sum_{j=0}^{k}\binom{k}{j}(-1)^{k-j}(j)_{n,\lambda},\quad (n\ge k\ge 0).
\end{displaymath} \par
Also, we have
\begin{align}
\sum_{n=k}^{\infty}S_{2,\lambda}(n,k)\frac{t^{n}}{n!}&=\frac{1}{k!}\big(e_{\lambda}(t)-1\big)^{k}=\frac{1}{k!}\bigg(\sum_{i=1}^{\infty}\frac{(1)_{i,\lambda}}{i!}t^{i}\bigg)^{k} \label{20} \\
&=\sum_{n=k}^{\infty}B_{n,k}\big((1)_{1,\lambda},(1)_{2,\lambda},\dots,(1)_{n-k+1,\lambda}\big)\frac{t^{n}}{n!}. \nonumber	
\end{align}
Comparing the coefficients on both sides of \eqref{20}, we have
\begin{equation}
S_{2,\lambda}(n,k)=B_{n,k}\big((1)_{1,\lambda},(1)_{2,\lambda},\dots,(1)_{n-k+1,\lambda}\big).\label{21}
\end{equation} \par
We note that
\begin{align}
\frac{1}{2}(e_{\lambda}(t)+1)&=\frac{1}{2}\bigg(\sum_{l=0}^{\infty}\frac{(1)_{l,\lambda}}{l!}t^{l}+1\bigg) \label{22}\\
&=1+\frac{1}{2}\sum_{l=1}^{\infty}(1)_{l,\lambda}\frac{t^{l}}{l!}, \nonumber	
\end{align}
and
\begin{align}
B_{n,k}\Big(\frac{(1)_{1,\lambda}}{2}, \frac{(1)_{2,\lambda}}{2},\dots, \frac{(1)_{n-k+1,\lambda}}{2}\Big)&=\frac{1}{2^{k}}B_{n,k}\big((1)_{1,\lambda},\dots,(1)_{n-k+1,\lambda}\big) \label{23} \\
&=\Big(\frac{1}{2}\Big)^{k}S_{2,\lambda}(n,k). \nonumber
\end{align} \par
From the definition of degenerate higher-order Euler numbers, we note that
\begin{align}
\sum_{n=0}^{\infty}\mathcal{E}_{n,\lambda}^{(\alpha)}\frac{t^{n}}{n!}&=\Big(\frac{1}{2}\big(e_{\lambda}(t)+1\big)\Big)^{-\alpha}=\Big(1+\sum_{i=1}^{\infty}\frac{1}{2}(1)_{i,\lambda}\frac{t^{i}}{i!}\Big)^{-\alpha} \label{24} \\
&=1+\sum_{k=1}^{\infty}(-\alpha)_{k}\Big(\frac{1}{2}\Big)^{k}\frac{1}{k!}\Big(\sum_{i=1}^{\infty}\frac{(1)_{i,\lambda}}{i!}t^{i}\Big)^{k}\nonumber\\
&=1+\sum_{k=1}^{\infty}(-\alpha)_{k}\Big(\frac{1}{2}\Big)^{k}\sum_{n=k}^{\infty}B_{n,k}\big((1)_{1,\lambda},(1)_{2,\lambda},\dots,(1)_{n-k+1,\lambda}\big)\frac{t^{n}}{n!}\nonumber \\
&=1+\sum_{n=1}^{\infty}\sum_{k=1}^{n}(-\alpha)_{k}\Big(\frac{1}{2}\Big)^{k}B_{n,k}\big((1)_{1,\lambda},(1)_{2,\lambda},\dots,(1)_{n-k+1,\lambda}\big)\frac{t^{n}}{n!}. \nonumber
\end{align}
Therefore, by comparing the coefficients on both sides of \eqref{24}, we obtain the following theorem.
\begin{theorem}
For $n\ge 1$, we have
\begin{align*}
\mathcal{E}_{n,\lambda}^{(\alpha)}&=\sum_{k=1}^{n}(-\alpha)_{k}\Big(\frac{1}{2}\Big)^{k}B_{n,k}\big((1)_{1,\lambda},(1)_{2,\lambda},\dots,(1)_{n-k+1,\lambda}\big) \\
&=\sum_{k=1}^{n}(-\alpha)_{k}\Big(\frac{1}{2}\Big)^{k}S_{2,\lambda}(n,k).
\end{align*}
\end{theorem}
From Theorem 5, we have
\begin{align}
A_{n,\lambda}^{(r,\alpha)}(x)&=(n)_{r}(x)_{n-r,\lambda}+(n)_{r}\sum_{k=0}^{n-r-1}\binom{n-r}{k}\mathcal{E}_{n-k-r,\lambda}^{(\alpha)}(x)_{k,\lambda} \label{25}\\
&=(n)_{r}(x)_{n-r,\lambda}+(n)_{r}\sum_{k=0}^{n-r-1}\binom{n-r}{k}(x)_{k,\lambda}\sum_{j=1}^{n-k-r}(-\alpha)_{j}\bigg(\frac{1}{2}\bigg)^{j}S_{2,\lambda}(n-k-r,j) \nonumber \\
&=(n)_{r}(x)_{n-r,\lambda}+(n)_{r}\sum_{k=0}^{n-r-1}\sum_{j=1}^{n-k-r}\binom{n-k}{k}(-\alpha)_{j}\bigg(\frac{1}{2}\bigg)^{j}S_{2,\lambda}(n-k-r,j)(x)_{k,\lambda},\nonumber
\end{align}
and
\begin{align}
A_{n,\lambda}^{(r,\alpha)}=(n)_{r}\sum_{j=1}^{n-r}(-\alpha)_{j}\bigg(\frac{1}{2}\bigg)^{j}S_{2,\lambda}(n-r,j),\label{26}
\end{align}
where $n,r\in\mathbb{Z}$ with $n>r\ge 0$. \par
Replacing $n$ by $n+r$ in \eqref{25}, we get
\begin{equation}
A_{n+r,\lambda}^{(r,\alpha)}(x)=(n+r)_{r}(x)_{n,\lambda}+(n+r)_{r}\sum_{k=0}^{n-1}\sum_{j=1}^{n-k}(-\alpha)_{j}\bigg(\frac{1}{2}\bigg)^{j}\binom{n+r-k}{k}S_{2,\lambda}(n-k,j)(x)_{k,\lambda},\quad(n \ge 1). \label{27}	
\end{equation}
Thus, we have
\begin{equation}
(x)_{n,\lambda}=\frac{1}{(n+r)_{r}}A_{n+r,\lambda}^{(r,\alpha)}(x)-\sum_{k=0}^{n-1}\sum_{j=1}^{n-k}(-\alpha)_{j}\Big(\frac{1}{2}\Big)^{j}\binom{n+r-k}{k}S_{2,\lambda}(n-k,j)(x)_{k,\lambda},\quad(n \ge 1). \label{28}
\end{equation}
Therefore, by \eqref{28}, we obtain the following theorem.
\begin{theorem}
For $n\ge 1$, we have
\begin{displaymath}
(x)_{n,\lambda}=\frac{1}{(n+r)_{r}}A_{n+r,\lambda}^{(r,\alpha)}-\sum_{k=0}^{n-1}\sum_{j=1}^{n-k}(-\alpha)_{j}\Big(\frac{1}{2}\Big)^{j}\binom{n+r-k}{k}S_{2,\lambda}(n-k,j)(x)_{k,\lambda}.
\end{displaymath}
\end{theorem}
For $m\in\mathbb{N}$ with $m\equiv 1\ (\mathrm{mod}\ 2)$, we have
\begin{align}
&\sum_{n=0}^{\infty}\sum_{k=0}^{m-1}(-1)^{k}A_{n,\lambda}^{(r,\alpha)}\Big(\frac{x+k}{m}\Big)\frac{t^{n}}{n!}=\sum_{k=0}^{m-1}(-1)^{k}\sum_{n=0}^{\infty}A_{n,\lambda}^{(r,\alpha)}\Big(\frac{x+k}{m}\Big)\frac{t^{n}}{n!}\label{29}\\
&=\sum_{k=0}^{m-1}(-1)^{k}t^{r}\Big(\frac{2}{e_{\lambda}(t)+1}\Big)^{\alpha}e_{\lambda}^{\frac{x+k}{m}}(t)
=\Big(\frac{2}{e_{\lambda}(t)+1}\Big)^{\alpha}t^{r}e_{\lambda}^{\frac{x}{m}}(t)\sum_{k=0}^{m-1}(-1)^{k}e_{\lambda}^{\frac{k}{m}}(t) \nonumber \\
&=\Big(\frac{2}{e_{\lambda}(t)+1}\Big)^{\alpha}m^{r}\Big(\frac{t}{m}\Big)^{r}e_{m\lambda}^{x}\Big(\frac{t}{m}\Big)\frac{1+e_{\lambda}(t)}{1+e_{\lambda}^{1/m}(t)}\nonumber \\
&=\Big(\frac{2}{e_{\lambda}(t)+1}\Big)^{\alpha-1}m^{r}\Big(\frac{t}{m}\Big)^{r}e_{m\lambda}^{x}\Big(\frac{t}{m}\Big)\frac{2}{1+e_{m\lambda}\big(\frac{t}{m}\big)} \nonumber \\
&=\sum_{j=0}^{\infty}\mathcal{E}_{j,\lambda}^{(\alpha-1)}\frac{t^{j}}{j!}\sum_{l=0}^{\infty}A_{l,m\lambda}^{(r)}(x)\frac{1}{m^{l-r}}\frac{t^{l}}{l!} \nonumber \\
&=\sum_{n=0}^{\infty}\sum_{l=0}^{n}\binom{n}{l}A_{l,m\lambda}^{(r)}(x)\frac{1}{m^{l-r}}\mathcal{E}_{n-l,\lambda}^{(\alpha-1)}\frac{t^{n}}{n!}. \nonumber
\end{align}
Therefore, by comparing the coefficients on both sides of \eqref{29}, we obtain the following theorem.
\begin{theorem}
For $m\in\mathbb{N}$ with $m\equiv 1\ (\mathrm{mod}\ 2)$, we have
\begin{displaymath}
\sum_{k=0}^{m-1}(-1)^{k}A_{n,\lambda}^{(r,\alpha)}\Big(\frac{x+k}{m}\Big)= \sum_{l=0}^{n}\binom{n}{l}A_{l,m\lambda}^{(r)}(x)\frac{1}{m^{l-r}}\mathcal{E}_{n-l,\lambda}^{(\alpha-1)}.
\end{displaymath}
\end{theorem}
Note that $\mathcal{E}_{0,\lambda}^{(0)}=1$, and $\mathcal{E}_{n,\lambda}^{(0)}=0,\ (n\in\mathbb{N})$. \par
Let us take $\alpha=1$ in \eqref{29}. Then we have
\begin{equation}
\sum_{k=0}^{m-1}(-1)^{k}A_{n,\lambda}^{(r)}\Big(\frac{x+k}{m}\Big)=A_{n,m\lambda}^{(r)}(x)\Big(\frac{1}{m}\Big)^{n-r},\label{30}	
\end{equation}
where $m\in\mathbb{N}$ with $m\equiv 1\ (\mathrm{mod}\ 2)$. \par
The alternating degenerate power sum of integers $T_{k,\lambda}(n)$ is defined by
\begin{equation}
T_{k,\lambda}(n)=\sum_{i=0}^{n}(-1)^{i}(i)_{k,\lambda},\quad (k\in\mathbb{N}\cup\{0\}). \label{31}
\end{equation}
\noindent\emph{Remark.} $\displaystyle T_{k}(n)=\sum_{i=0}^{n}(-1)^{i}i^{k}\displaystyle$ is defined in [11]. \par
Note that
\begin{align}
\sum_{i=0}^{n}(-1)^{i}e_{\lambda}^{i}(t)&=\sum_{k=0}^{\infty}\sum_{i=0}^{n}(-1)^{i}(i)_{k,\lambda}\frac{t^{k}}{k!} \label{32} \\
&=\sum_{k=0}^{\infty}T_{k,\lambda}(n)\frac{t^{k}}{k!}, \nonumber
\end{align}
 and
\begin{equation}
\sum_{i=0}^{n}(-1)^{i}e_{\lambda}^{i}(t)=\frac{1-(-1)^{n+1}e_{\lambda}^{n+1}(t)}{e_{\lambda}(t)+1}=\frac{1}{2}\sum_{k=0}^{\infty}\big(\mathcal{E}_{k,\lambda}-(-1)^{n+1}\mathcal{E}_{k,\lambda}(n+1)\big)\frac{t^{k}}{k!}.\label{33}	
\end{equation}
By \eqref{32} and \eqref{33}, we get
\begin{equation}
T_{k,\lambda}(n)=\frac{1}{2}\big(\mathcal{E}_{k,\lambda}+(-1)^{n}\mathcal{E}_{k,\lambda}(n+1)\big),\quad (n,k\ge 0).
\end{equation} \par
It is easy to show that
\begin{equation}
1-(-1)^{n+1}e_{\lambda}^{n+1}(t)=1-(-1)^{n+1}\sum_{l=0}^{\infty}\frac{(n+1)_{l,\lambda}}{l!}t^{l}=\Big(1+(-1)^{n}\Big)+(-1)^{n}\sum_{i=1}^{\infty}\frac{(n+1)_{i,\lambda}}{i!}t^{i}. \label{35}	
\end{equation}
For a fixed $n$, let $\alpha_{k,\lambda}$ be the sequence with
\begin{equation}
\alpha_{0,\lambda}=1+(-1)^{n},\quad \alpha_{k,\lambda}=(-1)^{n}(n+1)_{k,\lambda},\quad (k\ge 1).
\end{equation}
Then we have
\begin{equation}
1-(-1)^{n+1}e_{\lambda}^{n+1}(t)=\sum_{k=0}^{\infty}\alpha_{k,\lambda}\frac{t^{k}}{k!}. \label{37}
\end{equation} \par
Now, we observe that
\begin{align}
&\frac{1}{2}\sum_{n=0}^{\infty}\mathcal{E}_{n,\lambda}\frac{t^{n}}{n!}=\frac{1}{2}\bigg(\frac{e_{\lambda}(t)+1}{2}\bigg)^{-1}=\frac{1}{2}\bigg(1+\frac{1}{2}\sum_{i=1}^{\infty}\frac{(1)_{i,\lambda}}{i!}t^{i}\bigg)^{-1} \label{38}\\
&=\frac{1}{2}+\frac{1}{2}\sum_{k=1}^{\infty}(-1)_{k}\Big(\frac{1}{2}\Big)^{k}\frac{1}{k!}\bigg(\sum_{i=1}^{\infty}\frac{(1)_{i,\lambda}}{i!}t^{i}\bigg)^{k} \nonumber \\
&=\frac{1}{2}+\frac{1}{2}\sum_{k=1}^{\infty}(-1)^{k}k! \Big(\frac{1}{2}\Big)^{k}\sum_{n=k}^{\infty}B_{n,k}\big((1)_{1,\lambda},\dots,(1)_{n-k+1,\lambda}\big)\frac{t^{n}}{n!}\nonumber \\
&=\frac{1}{2}+\frac{1}{2}\sum_{n=1}^{\infty}\sum_{k=1}^{n}k!\Big(-\frac{1}{2}\Big)^{k}B_{n,k}\big((1)_{1,\lambda},\dots,(1)_{n-k+1,\lambda}\big)\frac{t^{n}}{n!} \nonumber \\
&=\frac{1}{2}+\frac{1}{2}\sum_{n=1}^{\infty}\sum_{k=1}^{n}(-1)^{k}k!\Big(\frac{1}{2}\Big)^{k}S_{2,\lambda}(n,k)\frac{t^{n}}{n!}. \nonumber
\end{align}
From \eqref{38}, we see that
\begin{equation}
\mathcal{E}_{0,\lambda}=1,\quad \mathcal{E}_{n,\lambda}=\sum_{k=1}^{n}(-1)^{k}k!\Big(\frac{1}{2}\Big)^{k}S_{2,\lambda}(n,k),\quad (n\ge 1). \label{40}
\end{equation} \par
From \eqref{32} and \eqref{37}, we note that
\begin{align}
&\sum_{k=0}^{\infty}T_{k,\lambda}(n)\frac{t^{k}}{k!}=\sum_{i=0}^{n}(-1)^{i}e_{\lambda}^{i}(t)=\frac{1-(-1)^{n+1}e_{\lambda}^{n+1}(t)}{1+e_{\lambda}(t)} \label{41} \\
&=\sum_{l=0}^{\infty}\alpha_{l,\lambda}\frac{t^{l}}{l!}\sum_{j=0}^{\infty}\frac{1}{2}\mathcal{E}_{j,\lambda}\frac{t^{j}}{j!}=\sum_{k=0}^{\infty}\frac{1}{2}\sum_{j=0}^{k}\binom{k}{j}\mathcal{E}_{j,\lambda}\alpha_{k-j,\lambda}\frac{t^{k}}{k!}. \nonumber
\end{align}
Using \eqref{41}, for $k \ge 1$, we have

\begin{align}
T_{k,\lambda}(n)&=\frac{1}{2}\sum_{j=0}^{k}\binom{k}{j}\mathcal{E}_{j,\lambda}\alpha_{k-j,\lambda}=\frac{1}{2}\mathcal{E}_{0,\lambda}\alpha_{k,\lambda}+\frac{1}{2}\mathcal{E}_{k,\lambda}\alpha_{0,\lambda}+\frac{1}{2}\sum_{j=1}^{k-1}\binom{k}{j}\mathcal{E}_{j,\lambda}\alpha_{k-j,\lambda}\label{42} \\
&=\frac{(-1)^{n}}{2}(n+1)_{k,\lambda}+\big(1+(-1)^{n}\big)\sum_{j=1}^{k}(-1)^{j}j!\bigg(\frac{1}{2}\bigg)^{j+1}S_{2,\lambda}(k,j)\nonumber \\
&\quad +\sum_{j=1}^{k-1}\binom{k}{j}\sum_{i=1}^{j}(-1)^{i}i!\bigg(\frac{1}{2}\bigg)^{i+1}S_{2,\lambda}(j,i)(-1)^{n}(n+1)_{k-j,\lambda} \nonumber \\
&=\frac{(-1)^{n}}{2}(n+1)_{k,\lambda}+\big(1+(-1)^{n}\big)\sum_{j=1}^{k}(-1)^{j}j!\bigg(\frac{1}{2}\bigg)^{j+1}S_{2,\lambda}(k,j) \nonumber \\
&\quad +(-1)^{n}\sum_{j=1}^{k-1}\sum_{i=1}^{j}\binom{k}{j}(-1)^{i}i!\bigg(\frac{1}{2}\bigg)^{i+1}S_{2,\lambda}(j,i)(n+1)_{k-j,\lambda}.\nonumber	
\end{align}

Therefore, by \eqref{42}, we obtain the following theorem.

\begin{theorem}
For $n \ge 0$ and $k \ge 1$, we have
\begin{align}
T_{k,\lambda}(n)&=\frac{(-1)^{n}}{2}(n+1)_{k,\lambda}+\big(1+(-1)^{n}\big)\sum_{j=1}^{k}(-1)^{j}j!\Big(\frac{1}{2}\Big)^{j+1}S_{2,\lambda}(k,j)\nonumber \\
&\quad +(-1)^{n}\sum_{j=1}^{k-1}\binom{k}{j}(n+1)_{k-j,\lambda}\sum_{i=1}^{j}(-1)^{i}i!\Big(\frac{1}{2}\Big)^{i+1}S_{2,\lambda}(j,i). \nonumber	
\end{align}
\end{theorem}
From Theorem 9, we note that
\begin{align}
T_{k,\lambda}(2n)&=\frac{1}{2}(2n+1)_{k,\lambda}+2\sum_{j=1}^{k}(-1)^{j}j!\Big(\frac{1}{2}\Big)^{j+1}S_{2,\lambda}(k,j) \label{43} \\
&\quad + \sum_{j=1}^{k-1}\binom{k}{j}(2n+1)_{k-j,\lambda}\sum_{i=1}^{j}(-1)^{i}i!
\Big(\frac{1}{2}\Big)^{i+1}S_{2,\lambda}(j,i),\nonumber
\end{align}
and
\begin{equation}
T_{k,\lambda}(2n+1)=-2^{k-1}(n+1)_{k,\lambda/2}-\sum_{j=1}^{k-1}\binom{k}{j}(n+1)_{k-j,\lambda/2}\sum_{i=1}^{j}(-1)^{i}i!2^{k-j-i-1}S_{2,\lambda}(j,i). \label{44}
\end{equation} \par
From \eqref{25}, we have
\begin{align}
A_{n,\lambda}^{(r)}(x)&=(n)_{r}(x)_{n-r,\lambda}+(n)_{r}\sum_{k=0}^{n-r-1}\sum_{j=1}^{n-k-r}\binom{n-r}{k}(-1)_{j}\Big(\frac{1}{2}\Big)^{j}S_{2,\lambda}(n-k-r,j)(x)_{k,\lambda} \label{45} \\
&=(n)_{r}	(x)_{n-r,\lambda}+(n)_{r}\sum_{k=0}^{n-r-1}\sum_{j=1}^{n-k-r}\binom{n-r}{k}j!\Big(-\frac{1}{2}\Big)^{j}S_{2,\lambda}(n-k-r,j)(x)_{k,\lambda}. \nonumber
\end{align}
Thus, by \eqref{45}, we get
\begin{equation}
A_{n,\lambda}^{(r)}\Big(\frac{x+i}{m}\Big)=(n)_{r}\Big(\frac{x+i}{m}\Big)_{n-r,\lambda}+(n)_{r}\sum_{k=0}^{n-r-1}\sum_{j=1}^{n-k-r}\binom{n-r}{k}j!\Big(-\frac{1}{2}\Big)^{j}S_{2,\lambda}(n-k-r,j)\Big(\frac{x+i}{m}\Big)_{k,\lambda}. \label{46}	
\end{equation}
From \eqref{46}, we can derive the following equation \eqref{47}.
\begin{align}
	\sum_{i=0}^{m-1}(-1)^{i}A_{m,\lambda}^{(r)}\Big(\frac{x+i}{m}\Big)&=(n)_{r}\sum_{i=0}^{m-1}(-1)^{i}\Big(\frac{x+i}{m}\Big)_{n-r,\lambda} \label{47} \\
	&\quad +(n)_{r}\sum_{k=0}^{n-r-1}\sum_{j=1}^{n-k-r}j!\Big(-\frac{1}{2}\Big)^{j}\binom{n-r}{k}S_{2,\lambda}(n-k-r,j)\sum_{i=0}^{m-1}(-1)^{i}\Big(\frac{x+i}{m}\Big)_{k,\lambda}. \nonumber
\end{align} \par
Now, we observe that
\begin{align}
\sum_{i=0}^{m-1}(-1)^{i}\Big(\frac{x+i}{m}\Big)_{k,\lambda}&=\frac{1}{m^{k}}\sum_{i=0}^{m-1}(-1)^{i}(x+i)_{k,m\lambda} \label{48}\\
&=\frac{1}{m^{k}}\sum_{i=0}^{m-1}(-1)^{i}\sum_{l=0}^{k}\binom{k}{l}(i)_{l,m\lambda}(x)_{k-l,m\lambda}\nonumber \\
&=\frac{1}{m^{k}}\sum_{l=0}^{k}\binom{k}{l}(x)_{k-l,m\lambda}T_{l,m\lambda}(m-1) \nonumber
\end{align}
From \eqref{47} and \eqref{48}, we have
\begin{align}
&\sum_{i=0}^{m-1}(-1)^{i}A_{n,\lambda}^{(r)}\Big(\frac{x+i}{m}\Big)=(n)_{r}\frac{1}{m^{n-r}}\sum_{l=0}^{n-r}\binom{n-r}{l}(x)_{n-l-r,m\lambda}T_{l,m\lambda}(m-1) \label{49}\\
&\quad +(n)_{r}\sum_{k=0}^{n-r-1}\sum_{j=1}^{n-k-r}j!\Big(-\frac{1}{2}\Big)^{j}	\binom{n-r}{k}S_{2,\lambda}(n-k-r,j)\frac{1}{m^{k}}\sum_{l=0}^{k}\binom{k}{l}T_{l,m\lambda}(m-1)(x)_{k-l,m\lambda}.\nonumber
\end{align} \par
For $m\in\mathbb{N}$ with $m\equiv 1\ (\mathrm{mod}\ 2)$, by Theorem 2 and \eqref{49}, we get
\begin{align}
A_{n,\lambda}^{(r)}(x)&=m^{n-r}\sum_{l=0}^{m-1}(-1)^{l}A_{n,\lambda/m}^{(r)}\Big(\frac{l+x}{m}\Big) \label{50}\\
&=(n)_{r}\sum_{l=0}^{n-r}\binom{n-r}{l}(x)_{n-l-r,\lambda}T_{l,\lambda}(m-1) \nonumber \\
&\quad +(n)_{r}\sum_{k=0}^{n-r-1}\sum_{j=1}^{n-k-r}\sum_{l=0}^{k}j!\Big(-\frac{1}{2}\Big)^{j}\binom{n-r}{k}\binom{k}{l}m^{n-r-k}S_{2,\lambda}(n-k-r,j)T_{l,\lambda}(m-1)(x)_{k-l,\lambda}. \nonumber
\end{align}
Therefore, by \eqref{50}, we obtain the following theorem.
\begin{theorem}
For $m\in\mathbb{N}$ with $m\equiv 1\ (\mathrm{mod}\ 2)$, we have
\end{theorem}
\begin{align}
A_{n,\lambda}^{(r)}(x)&=(n)_{r}\sum_{l=0}^{n-r}\binom{n-r}{l}(x)_{n-l-r,\lambda}T_{l,\lambda}(m-1) \nonumber \\
&\quad +(n)_{r}\sum_{k=0}^{n-r-1}\sum_{j=1}^{n-k-r}\sum_{l=0}^{k}j!\Big(-\frac{1}{2}\Big)^{j}\binom{n-r}{k}\binom{k}{l}m^{n-r-k}S_{2,\lambda}(n-k-r,j)T_{l,\lambda}(m-1)(x)_{k-l,\lambda}. \nonumber
\end{align}

\section{conclusion}
In recent years, various degenerate versions of many special numbers and polynomials have been explored by using different methods as aforementioned in the introduction. Many   \par
In this paper, we introduced the generalized degenerate Euler-Genocchi polynomials as a degenerate version of the Euler-Genocchi polynomials. In addition, we introduced their higher-order version, namely the generalized degenerate Euler-Genocchi polynomials of order $\alpha$, as a degenerate version of the generalized Euler-Genocchi polynomials of order $\alpha$. Then we studied certain properties and identities involving those polynomials, the generalized falling factorials, the degenerate Euler polynomials of order $\alpha$, the degenerate Stirling numbers of the second kind, and the alternating degenerate power sum of integers.  \par
It is one of our future projects to continue to study various degenerate versions of some special numbers and polynomials and to find their applications to physics, science and engineering as well as to mathematics.

\vspace{0.2in}
{ Funding}

This work was supported by the Basic Science Research Program, the National Research Foundation of Korea, (NRF-2021R1F1A1050151).


\begin{thebibliography}{9}
\bibitem{1}
Aydin, M. S.; Acikgoz, M.; Araci, S. \emph{A new construction on the degenerate Hurwitz-zeta function associated with certain applications.} Proc. Jangjeon Math. Soc. \textbf{25} (2022), no. 2, 195–-203.
\bibitem{2}
Belbachir, H.; Hadj-Brahim, S.; Rachidi, M. \emph{Another determinental approach for a family
of Appell polynomials.} Filomat \textbf{12} (2018), 4155-4164.
\bibitem{3}
Belbachir, H.; Hadj-Brahim, S. \emph{Some explicit formulas for Euler-Genocchi polynomials.}
Integers \textbf{19} (2019), $\sharp$A28, 14 pp.
\bibitem{4}
Carlitz, L. \emph{Degenerate Stirling, Bernoulli and Eulerian numbers.} Utilitas Math. \textbf{15} (1979), 51--88.
\bibitem{5}
Carlitz, L. \emph{A degenerate Staudt-Clausen theorem.} Arch. Math. (Basel) \textbf{7} (1956), 28--33.
\bibitem{6}
Comtet, L. \emph{Advanced combinatorics.} The art of finite and infinite expansions. Revised and enlarged edition. D. Reidel Publishing Co., Dordrecht, 1974. xi+343 pp. ISBN: 90-277-0441-4-05-02
\bibitem{7}
Goubi, M. \emph{On a generalized family of Euler-Genocchi polynomials.} Integers \textbf{21} (2021), Paper No. A48, 13 pp.
\bibitem{8}
 Kim, B. M.; Jang, L.-C.; Kim, W.; Kwon, H.-I. \emph{Degenerate Changhee-Genocchi numbers and polynomials.} J. Inequal. Appl. 2017, Paper No. 294, 10 pp.
\bibitem{9}
Kim, D. S.; Kim, T. \emph{A note on a new type of degenerate Bernoulli numbers.} Russ. J. Math. Phys. \textbf{27} (2020), no. 2, 227--235.
\bibitem{10}
Kim, D. S.; Kim, T.; Lee, S.-H.; Park, J.-W. \emph{Some new formulas of complete and incomplete degenerate Bell polynomials.} Adv. Difference Equ. 2021, Paper No. 326, 10 pp.
\bibitem{11}
Kim, T. \emph{On the alternating sums of powers of consecutive integers.} J. Anal. Comput. \textbf{1} (2005), no. 2, 117--120.
\bibitem{12}
Kim, T.; Kim, D. S. \emph{On some degenerate differential and degenerate difference operators.} Russ. J. Math. Phys. \textbf{29} (2022), no. 1, 37--46.
\bibitem{13}
Kim, T.; Kim, D. S. \emph{Some identities on truncated polynomials associated with degenerate Bell polynomials.} Russ. J. Math. Phys. \textbf{28} (2021), no. 3, 342--355.
\bibitem{14}
 Kim, T.; Kim, D. S.; Jang, G.-W. \emph{On degenerate central complete Bell polynomials.} Appl. Anal. Discrete Math. \textbf{13} (2019), no. 3, 805--818.
\bibitem{15}
Kim, T.; Kim, D. S.; Jang, L.-C.; Lee, H.; Kim, H.-Y. \emph{Complete and incomplete Bell polynomials associated with Lah-Bell numbers and polynomials.} Adv. Difference Equ. 2021, Paper No. 101, 12 pp.
\bibitem{16}
Kim, T.; Kim, D. S.; Kwon, J.; Lee, H. \emph{Some identities involving degenerate $r$-Stirling numbers.} Proc. Jangjeon Math. Soc. \textbf{25} (2022), no. 2, 245--252.
\bibitem{17}
Kim, T.; Kim, D. S.; Kwon, J.; Lee, H.; Park, S.-H. \emph{Some properties of degenerate complete and partial Bell polynomials.} Adv. Difference Equ. 2021, Paper No. 304, 12 pp.
\bibitem{18}
Kim, T.; Kim, D. S.; Kwon, J.; Park, S.-H. \emph{Representation by degenerate Genocchi polynomials.} J. Math. 2022, Art. ID 2339851, 11 pp.
\bibitem{19}
Kwon, H. I.; Jang, L.-C.; Kim, D. S.; Seo, J.-J. \emph{On modified degenerate Genocchi polynomials and numbers.} J. Comput. Anal. Appl. \textbf{23} (2017), no. 3, 521--529.
\bibitem{20}
Lee, D. S.; Kim, H. K. \emph{On the new type of degenerate poly-Genocchi numbers and polynomials.} Adv. Difference Equ. 2020, Paper No. 431, 15 pp.
\bibitem{21}
Lim, D. \emph{Some identities of degenerate Genocchi polynomials.} Bull. Korean Math. Soc. \textbf{53} (2016), no. 2, 569--579.
\bibitem{22}
Usman, T.; Aman, M.; Khan, O.; Nisar, K. S.; Araci, S. \emph{Construction of partially degenerate Laguerre-Genocchi polynomials with their applications.} AIMS Math. \textbf{5} (2020), no. 5, 4399--4411.
\end{thebibliography}
\end{document}